\documentclass{gen-j-l}

\usepackage{amssymb}

\newtheorem{theorem}{Theorem}[section]

\theoremstyle{definition}
\newtheorem{definition}[theorem]{Definition}

\theoremstyle{remark}

\numberwithin{equation}{section}

\newcommand{\ve}{\varepsilon}

\begin{document}

\title{Rough path theory and stochastic calculus}


%    author one information
\author{Yuzuru Inahama}
\address{
Graduate School of Mathematics, Kyushu University,
Motooka 744, Nishi-ku, Fukuoka, 819-0395, Japan
}
\curraddr{}
\email{inahama@math.kyushu-u.ac.jp
}
\thanks{}

%    author two information
%\author{}
%\address{}
%\curraddr{}
%\email{}
%\thanks{}

%    \subjclass is required by all journals except JAG.
\subjclass[2010]{Primary }

\date{}

\dedicatory{}

%    The "communicated by" line appears only in PROC and JAG.
%\commby{}

%    Abstract is required.
\begin{abstract}
T. Lyons' rough path theory is something like a deterministic version 
of K. It\^o's theory of  stochastic differential equations, 
combined with ideas from K. T. Chen's theory of iterated path integrals.
In this article we survey rough path theory,
in particular, its probabilistic aspects.
\end{abstract}

\maketitle


\bibliographystyle{amsplain}

%%%%%%%%%%%%%%%%%%%%%
\section{Introduction}

This article is a brief survey on rough path theory,
in particular, on its probabilistic aspects.
In the first half, 
we summarize basic results 
in the deterministic part of the theory.
The most important among them are
ODEs in rough path sense.
In the latter half, we discuss several important probabilistic results
in the theory.
Though putting them all in 
 a short article like this is not so easy, 
we believe it is worth trying because of importance and 
potential of rough path theory.

In 1998 T. Lyons \cite{ly} invented  rough path theory 
and then he wrote a book \cite{lq} with Z. Qian
which contains early results on rough paths.
This book is splendid mathematically.
However, because of minor errors 
and its very general setting, this book is not so readable.
Therefore, 
it was not easy to learn this theory
for non-experts who wanted to enter this research area.
(A few other books were published after that
and the situation has changed.
See Lyons, Caruana and Levy \cite{lcl}, Friz and Victoir \cite{fv},
Friz and Hairer \cite{fh}.)
Unlike in these thick standard books, 
in this article we will try to give a brief overview
of rough path theory
without  computations and proofs
so that the reader could grasp what the theory is all about.

A sample path of Brownian motion
is an important example of continuous paths 
in probability theory, 
but its behavior is quite bad.
In this theory, 
which  has one of its roots 
in K. T. Chen's theory of iterated path integrals, 
(objects corresponding to) 
iterated integrals of such bad paths are considered.
As a result, 
line integrals along a path or 
ordinary differential equations (ODEs) driven by a path
are generalized.
This, in turn, makes pathwise study of
stochastic differential equations (SDEs) possible.
In other words, T. Lyons successfully "de-randomized"
the SDE theory.
In particular, he proved that 
a solution to an SDE, 
as a functional of driving Brownian motion, 
becomes continuous.
From the viewpoint of the standard SDE theory
in which the martingale integration theory is crucially used,
this is quite surprising.
The SDE theory is very important and 
and 
has been a central topic  in probability theory
without exaggeration.
Since it has a long history and has been 
intensively and extensively studied 
by so many researchers,
this research area looked somewhat mature and 
some experts may have had a feeling that  
no big progress would be made
 when rough path theory was invented.

Rough path theory looks at SDEs from a very different angle
and we believe that 
it is breaking through the above-mentioned situation.
The number of researchers were not large, 
but 
it started to increase around 2010
as well as the number of papers.
In retrospect, 
this was probably when rough path theory really "took off."
Since this research area is still young, 
there will probably be many chances left for newcomers.
Indeed, we still saw unexpected developments recently,
which indicates that 
the theory is quite active and has large potential.   
%
%Since this research area is still young, 
%there will probably be many chances left for newcomers.
%
%
The purpose of this article is to give a bird's eye view 
of the rough path world to those who wants to enter it
and  to everyone who is interested 
in the theory, too.

%%%%%%%%%%%%%%%%%%%%%%%
\section{What is rough path theory?}

Though rough path theory is rapidly developing,
the population of researchers are not very large.
Even among probabilists,
not so many seem to
understand the outline of the theory.
Therefore, our aim of this section is to give a heuristic 
explanation on 
what the theory is all about.
The contents of this section are not intended to be rigorous
and  small matters are left aside.

Let us start with an ordinary differential equation (ODE)
driven by a path.
This type of ODE 
is usually called a driven ODE or a controlled ODE.
Let $x: [0,1] \to {\mathbf R}^d$
be a "sufficiently nice" path that starts at the origin.
(In this article all paths are continuous).
Let $\sigma: {\mathbf R}^n \to {\rm Mat} (n,d)$
and 
$b: {\mathbf R}^n \to {\mathbf R}^n$
be sufficiently nice functions,
where ${\rm Mat} (n,d)$ stands for the set of 
$n \times d$ real matrices.
Consider the following  ODE driven by the path $x$:
\[
dy_t = \sigma(y_t) dx_t +b (y_t)dt   \qquad \mbox{with given  $y_0 \in {\bf R}^n$.}
\]
This is slightly informal and its precise 
definition should be given by the following integral equation:
\[
y_t = y_0 + \int_0^t \sigma(y_s) dx_s +  \int_0^t b (y_s)ds .
\]
When there exists a unique solution,
$y$ can be regarded as a function (or a map)
of $x$.
Using the terminology of probability theory, 
we call it the It\^o map.
It is a map from one path space to another.
We will assume for simplicity 
that $y_0 =0$ and $b \equiv 0$ 
(because we can take the space-time path 
$t \mapsto (x_t, t)$ and a block matrix 
$[\sigma | b]$ of size $n \times (d+1)$).
Therefore, 
we will consider 
\begin{equation}\label{eq.xode1}
dy_t = \sigma(y_t) dx_t \quad \mbox{with $y_0 =0$}
\qquad
\Longleftrightarrow
\qquad
y_t =  \int_0^t \sigma(y_s) dx_s 
\end{equation}
from now on.

Whether ODE (\ref{eq.xode1}) makes sense or not 
depends on 
well-definedness of the line integral on the right hand side.
If it is well-defined, 
then under a suitable condition on the regularity of 
the coefficient matrix $\sigma$,
we can usually obtain a (time-)local unique 
solution.
The most typical method is Picard's iteration 
on a shrunk time interval.

Note that for a generic continuous path,
the line integral cannot be defined.
A stronger condition on $x$ is needed.
For instance, for a piecewise $C^1$ path $x$,
the line integral clearly makes sense since 
$dx_s = x^{\prime}_s ds$.
A more advanced example could be 
a path of bounded variation.
In this case, the integral can be understood in
the Riemann-Stieltjes sense
and 
ODE (\ref{eq.xode1}) has a unique solution.
Moreover, 
if the path spaces are equipped with the bounded
variation norm, 
then  line integrals and It\^o maps become 
continuous as maps between path spaces.
These are basically within 
an advanced course of calculus and not so difficult.

Less widely known is the Young integral, 
which is essentially 
a generalized Riemann-Stieltjes integral.
We will briefly explain it below.
If $x$ is of finite $p$-variation and
$y$ is of finite $q$-variation 
with $p, q \ge 1$ and $1/p +1/q >1$,
then the line integral in   (\ref{eq.xode1}) makes sense.
The approximating Riemann sum that defines the Young integral
is 
exactly the same as the one for the Riemann-Stieltjes integral.
It is obvious from this that 
the Young integral extends 
the Riemann-Stieltjes integral if it exists.
We often use the Young integration theory with 
$p=q$.
In such a case, the Young integral 
is well-defined if $1 \le p=q <2$ and
with respect to $p$-variation norm ($1 \le p <2$),
 line integrals and It\^o maps become 
continuous, too.

However, for some reason that will be explained shortly, 
it cannot be used for 
stochastic integrals along Brownian motion.
Our main interest in this article is that 
"How far can we extend line integrals 
beyond Young's theory in a deterministic way
so that it can be used for probabilistic studies."
Of course, the main example we have in mind 
is a sample path of Brownian motion.

Denote by $\mu$ the $d$-dimensional Wiener measure,
that is,
the law of $d$-dimensional Brownian motion.
It sits on the space of continuous functions
$C_0 ({\bf R}^d)=\{ x: [0,1] \to {\bf R}^d ~|~ \mbox{conti, $x_0=0$}\}$
and is the most important probability measure 
in probability theory.
The path $t \mapsto x_t$ can be viewed as 
a random motion under $\mu$.
In that case we call $(x_t)_{t \ge 0}$
(the canonical realization of) Brownian motion.
It is well-known that 
Brownian motion is a very zig-zag movement 
and its trajectory is very wild.
For example, for any $p \le 2$,
the set of paths with finite $p$-variation is 
a $\mu$-zero set.
\footnote{When standard textbooks on probability say that 
the quadratic variation of one-dimensional
 Brownian motion on $[0,T]$ equals $T$, 
the definition of "quadratic variation" is different from 
the one of the $2$-variation norm in this article.
So, there is no contradiction.}
Therefore, it is impossible to define line integral 
along Brownian paths by using 
the Young (or the Riemann-Stieltjes) integral.

In standard probability theory,
a line integral along Brownian paths is defined 
as It\^o's stochastic integral as follows 
(for simplicity we set $d=1$):
\[
 \int_0^t  z_s dx_s 
 =
 \lim_{|{\mathcal P}| \to 0}
 \sum_{i=1}^N  z_{t_{i-1}}  (x_{t_i} -x_{t_{i-1}}),
  \]
where 
${\mathcal P} =\{ 0=t_0 <t_1 <\cdots < t_N =t \}$
is a partition of $[0,1]$.
In defining 
and proving basic properties of this stochastic integral,
the martingale property of Brownian motion plays 
a crucial role.
Using it one can show that
\[
{\mathbf E} [ \bigl| \int_0^t z_s dx_s \bigr|^2 ]  = {\mathbf E} [  \int_0^t  |z_s|^2 ds ],
\]
which means that stochastic integration 
is 
an isometry 
between $L^2 (\mu \times ds)$ and $L^2 (\mu)$.
This is the most important fact 
in It\^o's theory of stochastic integration.

Each element of 
$L^2 (\mu)$ is just an equivalence class 
with respect to $\mu$
and a single-point set is of $\mu$-zero set, 
the stochastic integral 
does not have an $x$-wise meaning.
Neither is it continuous in $x$.
For example, 
let us consider L\'evy's stochastic area for two-dimensional 
Brownian motion
\begin{equation}\label{levy}
x=(x^1,x^2) \mapsto
\int_0^1  (  x^2_s dx^1_s -   x^1_s dx^2_s ).
\end{equation}
With respect to any Banach which preserves 
the Gaussian structure of the classical Wiener space,
the above map is discontinuous
(see Sugita \cite{su}).
As a result, the solution $y$ to equation (\ref{eq.xode1})
understood in the It\^o sense
is not continuous in $x$.
In other words, 
the It\^o map is not or cannot be made continuous 
in the driving path.

As we have seen, 
deterministic line integrals such as  the Young integral
have a limit
and are unsatisfactory from a probabilistic view point.
On the other hand, 
the It\^o integral turned out to be extremely successful.
In such a situation, 
discontinuity of the It\^o map and
impossibility of a pathwise definition of stochastic integrals 
were
probably "unpleasant facts one has to accept"
for most of the probabilists. 
This was the atmosphere in the probability community.

T. Lyons \cite{ly} made a breakthrough 
 by inventing rough path theory.  
It enables us to do pathwise study of SDEs.
In fact, in this theory
we consider not just a path itself, 
but also iterated integrals of the paths together.
A generalized path in this sense is called a rough path.
This idea probably comes from K. T. Chen's theory
of iterated integrals of paths in topology.
Unlike in topology, however,
we have to deal with paths with low regularity,
since our main interest is in probabilistic applications.
Therefore, we have to take completion of 
the set of nice paths with respect to a certain 
 Banach norm,
 but it is difficult to find a suitable norm.

 The most important feature of the theory is as follows:
"If  the rough path space is equipped with a 
suitable topology, 
then line integrals and It\^o maps can be defined 
in a deterministic way 
and they become continuous." 
The continuity of It\^o maps is called Lyons' continuity 
theorem (or the universal limit theorem)
and is the pivot of the theory.

As is mentioned above, 
a rough path is a pair of 
its first and  second level paths.
The first level path is just a difference 
of a usual path (which is a single line integral of the path)
and 
the second level path is a double integral of the usual path.
If we agree that 
the starting point of a usual path is always the origin, 
then the path itself and its difference is equivalent.
So, a first level path is actually a path in the usual sense.
Novelty is in taking iterated integrals of a path
into consideration.

Choose  $2 <p <3$
and introduce a topology on the rough path space 
so that 
the first level paths are of finite $p$-variation
and 
the second level paths are of finite $p /2$-variation.
Then, it is important that
the following two seemingly opposite requests
are satisfied simultaneously. 
{\bf (a)}~Line integrals along a rough path 
can be defined deterministically.
This means that 
regularity of rough paths is nice 
and hence the rough path space is small
at least to this extent.
{\bf (b)}~(Lift of) Wiener measure sits on the rough path space.
This means 
that  the rough path space is large at least to this extent.

If we substitute the lift of Brownian motion 
in the Lyons-It\^o map, which is a rough path version 
of the It\^o map,
then we obtain the solution to the corresponding 
SDE of Stratonovich type.
(An SDE of  Stratonovich type is a slight modification of
an SDE of It\^o type.)
Recall that the driven ODE in rough path sense
is deterministic and irrelevant to any measure.
Therefore, SDEs are "de-randomized."
In other words, 
probability measures and driven ODEs are separated.
This is impossible as long as we use the martingale 
integration theory.

Before we end this section, 
we make clear what are basically 
not used in rough path theory.
{\bf (a)}~Martingale integration theory,
{\bf (b)}~Markov property,
{\bf (c)}~filtration, which is an increasing family 
of sub-$\sigma$ field indexed by the time parameter.
Consequently, 
this theory has a strong taste of real analysis and
 does not look like probability theory very much.

%%%%%%%%%%%%%%%%%%

\section{Geometric rough path}

In this section 
we define rough paths,  following \cite{lq, lcl, fvbk}. 
For simplicity 
we consider the case $2 \le p <3$,
where $p$ is a constant called roughness
and stands for the index of variation norm.
This is enough 
for applications to Brownian motion.
In this case only the first and the second level paths appear.
Some people prefer 
$1/p$-H\"older norm, which is a twin sister 
of $p$-variation norm,
but we basically use the variation norm in this article. 
Of course, 
rough path theory extends to the case $p \ge 3$.
In that case 
paths up to $[p]$th level,
which roughly correspond to 
$i$th iterated integral ($1 \le i \le [p]$) of the first level path,
 are used.

Set 
$\triangle:=\{ (s,t) ~|~ 0 \le s \le t \le 1 \}$.
For $p \ge 1$
and a continuous map
$A:\triangle \to {\bf R}^d$,
we define $p$-variation norm of $A$ by
\begin{equation}\label{def.norm}
\|A\|_{p} :=
\sup_{ {\mathcal P}}  \Bigl\{ \sum_{i} | A_{t_{i-1},t_i}|^p
\Bigr\}^{1/p}.
\end{equation}
Here, the supremum runs over all the finite partitions
${\mathcal P}=\{0 =t_0 <t_1 < \cdots <t_N =1 \}$
of 
$[0,1]$.
Note that,  if $p<p'$, then $\|A\|_{p} <\infty$ implies 
$\|A\|_{p'} <\infty$.
In other words, the larger $p$ is, the weaker the condition 
of finite $p$-variation becomes.
In particular,  if something is of finite $1$-variation,
it is considered to be "very nice" in this theory.
If you prefer the H\"older norm,
then instead of (\ref{def.norm})
use $\|A\|_{1/p -Hld} :=\sup_{ s<t} | A_{s,t}|/|t-s|^{1/p}$.

Let 
$T^{(2)}( {\bf R}^d) :=
{\bf R} \oplus {\bf R}^d\oplus ({\bf R}^d \otimes {\bf R}^d)$
be the truncated tensor algebra of degree $2$.
Now we define an ${\bf R}^d$-valued rough path
of roughness $p$.
The totality of such rough paths will be denoted
by $\Omega_p ({\bf R}^d)$.

\begin{definition}
A continuous map 
$X=(1, X^1, X^2): \triangle\to T^{(2)} ( {\bf R}^d)$
is said to be a rough path 
if the following two conditions are satisfied:
\\
\noindent
{\bf (i)}~(Chen's identity)~  
For any $0\le s \le u \le t \le 1$,
\[
X^1_{s,t} = X^1_{s,u}+X^1_{u,t},
\qquad
X^2_{s,t} = X^2_{s,u}+X^2_{u,t}+ X^1_{s,u}\otimes X^1_{u,t}.
\]
\noindent
{\bf (ii)}~(finite $p$-variation)
\qquad
$\|X^1\|_{p} <\infty, \quad \|X^2\|_{p/2} <\infty$.
\end{definition}

We will basically omit the obvious $0$th component "$1$"
and simply write $X= (X^1, X^2)$.
The two norms in condition {\bf (ii)} naturally 
defines a distance on $\Omega_p ( {\bf R}^d)$
and makes it a complete metric space (but not separable).
The first level path $X^1$ is just a difference 
of the usual path in ${\bf R}^d$ with finite $p$-variation.
At first sight, 
Chen's identity for 
the second level path $X^2$ may look strange.
As we will see, however,
 $X^2$ is an abstraction
of the two-fold iterated integral of a nice usual path 
in ${\bf R}^d$.
If the multiplication of $T^{(2)} ( {\bf R}^d)$
is denoted by $\otimes$,
then Chen's identity reads $X_{s,t} = X_{s,u} \otimes X_{u,t}$.
(This is the relation for differences 
of a group-valued path.)

Now we give a natural example of rough path.
It is very important both theoretically
and practically.
For a continuous path 
$x:[0,1] \to {\bf R}^d$ of finite $1$-variation
that starts from $0$
and $(s,t) \in \triangle$, 
set  
\begin{eqnarray*}
X^1_{s,t} &=& \int_s^t dx_{t_1} = x_t -x_s, 
\nonumber
\\
X^2_{s,t} &=& 
\int_{s \le t_1 \le t_2 \le t}   dx_{t_1}  \otimes dx_{t_2}
=\int_s^t (x_u -x_s) \otimes dx_u,
%\qquad
%(s,t) \in \triangle
\end{eqnarray*}
Then, it is straight forward to check that
$X \in \Omega_p ({\bf R}^d)$.
It is called a smooth rough path above $x$
(or the natural lift of $x$).
Note that the Riemann-Stieltjes (or Young) integral is 
used to define $X^2$.
Hence, a generic continuous path cannot be lifted
in this way.

Since $\Omega_p ({\bf R}^d)$ is a bit too large,
we introduce the geometric rough path space.
This is the main path space in rough path theory
and plays a role of the classical Wiener space in 
usual probability theory.

\begin{definition}
A rough path that can be approximated by smooth rough paths
is called a geometric rough path.
The set of geometric rough paths is denoted by 
$G\Omega_p ({\bf R}^d )$, namely,
$G\Omega_p ( {\bf R}^d)
=
\overline{ \{ \mbox{${\bf R}^d$-valued smooth rough paths} \} }^{d_p} 
\subset \Omega_p ({\bf R}^d)$.
\end{definition}

By way of construction 
$G\Omega_p ({\bf R}^d )$ becomes a complete separable 
metric space.
There exist $X, Y \in G\Omega_p ({\bf R}^d )$
such that $X^1 =Y^1$, but $X^2 \neq Y^2$.
This means that the second level paths 
do have new information.
For $X \in G\Omega_p ({\bf R}^d )$,
the symmetric part of $X^2_{s,t}$ 
 is determined by the first level path
 since it is given by $(X^1_{s,t} \otimes X^1_{s,t})/2$.
Hence, all information of $X$ is contained in
$X^1$ and the anti-symmetric part of $X^2$.
The latter is also called L\'evy area
and has a similar form to (\ref{levy}).
Therefore, 
things like L\'evy area are built in the structure
of 
$G\Omega_p ({\bf R}^d )$ and 
continuity of L\'evy area as functions 
on $G\Omega_p ({\bf R}^d )$ is almost obvious.

For 
$X, Y \in G\Omega_p ({\bf R}^d)$
the addition "$X+Y$" cannot be defined in general.
However, a natural scalar action called the dilation exists.
Similarly, 
for $X \in G\Omega_p ({\bf R}^d)$
and $Y  \in G\Omega_p ({\bf R}^r)$,
a paired rough path $"(X,Y)" \in G\Omega_p ( {\bf R}^d \oplus  {\bf R}^r)$
cannot be defined in general, either.
However, 
if one of $X$ and $Y$ is a smooth rough path,
then both $X+Y$ and $(X, Y)$ 
can be defined naturally
since the "cross integrals" of $X$ and $Y$ 
is well-defined as Riemann-Stieltjes integrals.
(This paragraph is actually important).

In the definition of  geometric rough paths,
paths of finite $1$-variation 
and the Riemann-Stieltjes integral are used.
However, 
even if they are replaced by 
paths of finite $q$-variation  with $1 \le q <2$
and the Young integral, respectively,
the definition remains  equivalent. 
Similarly, the addition $X+Y$ and the pair $(X, Y)$
are in fact well-defined 
if one of $X$ and $Y$ are 
of finite $q$-variation  with $1 \le q <2$ and $1/p +1/q >1$.
Hence, 
$X+Y$ and $(X, Y)$ are called the Young translation (shift)
and the Young pairing, respectively.

Before closing this section, 
we give a sketch of higher level geometric rough paths.
Simply put, 
basically everything in this section still holds 
with possible minor modifications
when the roughness $p \ge 3$.
We have to modify the following points.
The truncated tensor algebra 
$T^{([p])}( {\bf R}^d)$ of degree $[p]$ is used.
The $i$th level path is estimated by $p/i$-variation norm 
($1 \le i \le [p]$).
When we lift a usual path $x$ of finite variation, 
we consider
\[
X^i_{s,t} =
\int_{s \le t_1 \le \cdots \le  t_i \le t}   dx_{t_1}  \otimes \cdots \otimes dx_{t_i}
\qquad
(1 \le i \le [p], \, (s,t)\in \triangle ),
\]
that is, 
all  iterated integrals of $x$ of degree up to $[p]$.
Chen's identity can be written 
as $X_{s,t} = X_{s,u} \otimes X_{u,t}$ as before,
which is the algebraic relation of 
differences of a group-valued path.
What is the smallest group which contain 
all such $X_{s,t}$'s as $x$ and $s,t$ vary?
The answer is 
the free nilpotent Lie group $G^{[p]}$ of step $[p]$,
which is a subgroup of $T^{([p])}( {\bf R}^d)$.
This group has a nice homogeneous distance
which is compatible with the dilation.
Once one understands basic properties of $G^{[p]}$ 
and this distance,
one can clearly see 
why the $i$th level path is estimated 
by the $p/i$-variation norm.
Loosely speaking, a geometric rough path 
is equivalent 
to a continuous path on $G^{[p]}$ staring at the unit
with finite $p$-variation with respect to this distance.
Therefore, a geometric rough path 
is never a bad object despite its looks.
This point of view is quite useful when $p$ is large.
(The contents of this paragraph is well summarized in
Friz and Victoir \cite{fvbk}.)
We remark that the geometric rough path space 
with $1/p$-H\"older topology is defined in a similar way.

%%%%%%%%%%%%%%%%%%%%%%%%%%
\section{Line integral along  rough path}
In this section 
we discuss line integrals along a rough path
when $2 \le p <3$. 
Let 
$X \in G\Omega_p ( {\bf R}^d)$
and
$f:{\bf R}^d \to {\rm Mat} (n,d)$ be of $C^3$
which should be viewed as a vector-valued $1$-form.
We would like to define an integral
$\int f(X)dX$
as an element of $G\Omega_p ( {\bf R}^n)$.
(The condition of $f$ can be relaxed slightly.)
Note that $\int f(X)dY$ cannot be defined in general
except when $(X,Y)$ defines a rough path over 
the direct sum space.
The contents of this section naturally 
extends to the case $p \ge 3$, too.

Now we introduce a Riemann sum which 
approximates the rough path integral.
We write $x_s =X^1_{0,s}$.  For $(s,t) \in \triangle$, we set
\begin{align}
\hat{Y}^1_{s,t}  &=  f( x_s ) X^1_{s,t}+  \nabla f( x_s ) X^2_{s,t},
%\label{def_int1}
\nonumber\\
\hat{Y}^2_{s,t}  &=  f( x_s )\otimes f( x_s ) X^2_{s,t}.
%\label{def_int2}
\nonumber
\end{align}
Here,
$\hat{Y}^1_{s,t} \in {\bf R}^n$ and $\hat{Y}^2_{s,t} \in {\bf R}^n \otimes {\bf R}^n$.
Note that if the second term on the right hand side were absent,
the first one would be
 just a summand for the usual Riemann sum.

Let 
${\mathcal P}=\{ s=t_0 <t_1< \cdots < t_n=t \}$
be a partition of $[s,t]$ and denote by
$|{\mathcal P}|$ its mesh.
If we set
\begin{eqnarray}
Y^1_{s,t}  &=&  \lim_{|{\mathcal P}| \searrow 0}
\sum_{i=1}^{n}  \hat{Y}^1_{t_{i-1},  t_{i}},
\label{def_int3}
\\
Y^2_{s,t}  &=&  \lim_{|{\mathcal P}| \searrow 0}
\sum_{i=1}^{n}
\bigl(
\hat{Y}^2_{t_{i-1},  t_{i}}  
+   Y^1_{s, t_{i-1}} \otimes Y^1_{t_{i-1},  t_{i}}
\bigr),
\label{def_int4}
\end{eqnarray}
then the right hand sides of the both equations
converge and it holds that
$Y=(Y^1, Y^2) \in G\Omega_p ({\bf R}^n)$.
\footnote{In fact, $\hat{Y}$ is an almost rough path is the sense of 
Lyons and Qian \cite{lq}.
For every almost rough path, there exists a unique rough path 
associated with it.
Equations (\ref{def_int3})--(\ref{def_int4}) are actually a special case 
of this general theorem.}
We usually write $Y^j_{s,t}= \int_{s}^{t} f(X)dX^j~(j=1,2)$.
At first sight  (\ref{def_int4}) may look strange,
but it is not.
To see this, one should first rewrite Chen's identity for 
not just two subintervals of $[s,t]$,
but for $n$ subintervals 
and then compare it to  (\ref{def_int4}).
With respect to the natural distances on the geometric 
rough path spaces,
the map $X \mapsto \int f(X)dX$ is 
locally Lipschitz continuous,
that is, Lipschitz continuous on any bounded set.

Let us summarize. 
\begin{theorem}
If $f:{\bf R}^d \to {\rm Mat}(n,d)$ is of $C^3$,
the rough path integration map
\[
G\Omega_p ( {\bf R}^d)  \ni X  \mapsto
\int f(X)dX  \in G\Omega_p ({\bf R}^n)
\]
is locally Lipschitz continuous
and extends
the Riemann-Stieltjes integration map
$x \mapsto \int_0^{\cdot} f(x_s) dx_s$.
\end{theorem}

Before ending this section,
we make a simple remark on  the rough path integration.
Without loss of generality  we assume $n=1$. 
Hence, $f$ is a usual one-form on ${\bf R}^d$.
It is obvious that if $f$ is exact, 
that is, $f = dg$ for a some function $g: {\bf R}^d \to {\bf R}$,
then $\int_0^{T} f(x_s) dx_s = g(x_T) -g(x_0)$
and the line integral clearly extends to any continuous path $x$.
Therefore, 
when one tries to extend  line integration,
non-exact one-forms are troubles.
The simplest non-exact one-forms on  ${\bf R}^d$ are 
$\xi_i d\xi_j - \xi_j d\xi_i ~(i<j)$, 
where $(\xi_1, \ldots, \xi_d)$ is the coordinate of ${\bf R}^d$.
The line integrals along a path $x$ of those one-forms
are the L\'evy areas of $x$.
Remember that information of the L\'evy areas is
 precisely what is added to a path 
when it gets lifted to a geometric rough path.
Therefore, the point to observe is this:
Even though only line integrals along $x$ of $\xi_i d\xi_j - \xi_j d\xi_i ~(i<j)$
were added,
 line integrals along $x$ of every one-form $f$ are 
 continuously extended.

%%%%%%%%%%%%%%%%%%%%%%%%%%%%
\section{ODE driven by  rough path}
In this section 
we consider a driven ODE in the sense of 
rough path theory 
(rough differential equation, RDE).
We follow Lyons and Qian \cite{lq}.
For simplicity, we assume $2 \le p <3$.
However, the results in this section holds for $p \ge 3$.
One should note that 
an RDE is deterministic.
In this section
$\sigma:  {\bf R}^n \to {\rm Mat } (n,d)$
is assume to be of $C^3_b$,
that is, $| \nabla^j \sigma |$ is bounded for $0 \le j \le 3$.

For a given ${\bf R}^d$-valued path $X$, 
we consider the following (formal) driven ODE:
\begin{equation}
dY_t = \sigma (Y_t) dX_t, \qquad Y_{0}=0.
\label{eq.ode1}
\end{equation}
A solution $Y$ is an ${\bf R}^n$-valued path.
When we consider a non-zero initial condition 
$Y_0=y_0 \in {\bf R}^n$,
we replace the coefficient by $\sigma ( \cdot +y_0)$.
As always this ODE should be defined as an integral equation:
\[
Y_t =\int_0^t  \sigma (Y_u) dX_u.
\]
In rough path theory,
however, the right hand side does not make sense
since $X$ and $Y$ are different rough paths 
and the rough path integral may be ill-defined. 
So we add a trivial equation to (\ref{eq.ode1})
and consider
the following system of ODEs instead:
\begin{eqnarray}
\left\{
\begin{array}{@{\,}lll}
dX_t &=& dX_t,  \\
dY_t &=&  \sigma (Y_t) dX_t.
\end{array}
\right.
%
%
%\nonumber
\label{eq.ode3}
\end{eqnarray}

The natural projections 
from a direct sum ${\bf R}^d \oplus  {\bf R}^n$ to 
each component is denoted by $\pi_1, \pi_2$,
respectively. 
Namely, $\pi_1 z =x$ and $\pi_2 z =y$
for $z=(x,y)$.
Define $\hat{\sigma}:{\bf R}^d \oplus {\bf R}^n
\to {\rm Mat}(d+n , d+n)$
by
\begin{eqnarray}
\hat{\sigma}(z)
=
\begin{pmatrix}
1 & 0 \\
\sigma( \pi_2 z) & 0  
\end{pmatrix}
\quad
\mbox{or }
\quad
\hat{\sigma}(z) \langle z' \rangle
=
\begin{pmatrix}
1 & 0 \\
\sigma( y) & 0 
\end{pmatrix}
\begin{pmatrix}
x' \\
y' 
\end{pmatrix}
=
\begin{pmatrix}
x' \\
\sigma( y)x'   
\end{pmatrix}.
\nonumber
%\label{eq.ode4}
\end{eqnarray}
Then,  (\ref{eq.ode3}) is equivalent to
\begin{eqnarray}
dZ_t = \hat{\sigma}(Z_t)dZ_t  \qquad \mbox{  with $\pi_1Z_t =X_t$.} 
%\label{eq.ode5}
\nonumber
\end{eqnarray}

Summarizing these, 
we set the following definition
(the initial value $y_0 =0$ is assumed).
The projection $\pi_1$ (resp. $\pi_2$) 
naturally induces a projection 
$G\Omega_p  ({\bf R}^d  \oplus {\bf R}^n) \to G\Omega_p ({\bf R}^d)$
(resp. $\to G\Omega_p ({\bf R}^n)$),
which will be denoted by the same symbol.

\begin{definition}\label{def_ode}
Let $X \in G\Omega_p ({\bf R}^d)$.
A geometric rough path 
$Z \in G\Omega_p ({\bf R}^d \oplus {\bf R}^n)$
is said to be a solution to (\ref{eq.ode1})
in the rough path sense if the following rough integral equation
is satisfied:
\begin{eqnarray}
Z = \int \hat{\sigma}(Z)dZ ,  \qquad\qquad \mbox{ with } \qquad  \pi_1Z =X 
\label{eq.ode6}
\end{eqnarray}
Note that
the second level $Y =\pi_2 Z$ is also called a solution.
If there is a unique solution, 
the map  $X \mapsto Y$ is called the Lyons-It\^o map
and denoted by $Y =\Phi (X)$. 
\end{definition}

As we have seen, 
in the original formalism of Lyons
a solution $Y$ does not exist alone, 
but is the second component of 
a solution rough path over a direct sum space.
In some new methods, however,
a solution to RDE is not defined to be  of the form $Z =(X,Y)$.
One of them is Gubinelli's formalism \cite{gu},
in which the driving rough path $X$ and 
a solution $Y$ are separated in a certain sense.

Now, we present the most important theorem 
in the theory, namely, Lyons' continuity theorem 
(also known as the universal limit theorem).
\begin{theorem}\label{thm_cont}
Consider RDE (\ref{eq.ode1})
with a $C^3_b$-coefficient 
$\sigma: {\bf R}^n \to {\rm Mat}(n,d)$.
Then, for any $X \in G\Omega_p ({\bf R}^d)$,
there exists a unique solution 
$Z \in G\Omega_p ( {\bf R}^d \oplus  {\bf R}^n)$
 to (\ref{eq.ode6}).
Moreover, $X \mapsto Z$ is locally Lipschitz continuous 
and so is the Lyons-It\^o map $X  \mapsto Y= \pi_2 Z=\Phi(X) \in G\Omega_p ( {\bf R}^n)$.
\end{theorem}

If $X$ is a smooth rough path lying above 
a usual ${\bf R}^d$-valued path $x$ with finite variation,
then $Y$ is a smooth rough path lying above 
a unique solution $y$ to the corresponding ODE 
in the Riemann-Stieltjes sense.
This can be easily
shown by the uniqueness of RDE in Theorem \ref{thm_cont}
and the fact that the rough path integration extends 
the Riemann-Stieltjes one. 
Thus, we have generalized driven ODEs.

\bigskip
\noindent
{\it Sketch of Proof of Theorem \ref{thm_cont}.}~
We use Picard's  iteration method.
Set $Z(0)$
by
$Z(0)^1_{s,t}=(X^1_{s,t},0), ~Z(0)^2_{s,t}=(X^2_{s,t},0,0,0)$
and set
\begin{eqnarray}
Z(m) =
\int \hat{\sigma}(Z(m-1))dZ(m-1)  
%
%\label{eq.ode6.5}
\nonumber
\end{eqnarray}
for $m \ge 1$.

If $T_1 \in (0,1]$ is small enough,
then the Lipschitz constant of 
the rough integration map in (\ref{eq.ode6})
becomes smaller than $1$.
Hence,  $\{ Z(m)\}_{m=0,1,2,\ldots}$
converges to some
$Z$ in 
$G\Omega_p ( {\bf R}^d \oplus  {\bf R}^n)$.
Thus, we find a solution on the subinterval $[0, T_1]$.

Next, we solve the RDE on $[T_1, T_2]$
with a new initial condition $y_{T_1}=Y^1_{0,T_1}$.
Repeating this, 
we obtain a solution on each subinterval $[T_{i},T_{i+1}]$.
This procedure stops for some finite $i$,
that is, $T_{i+1} \ge 1$,
because of the $C^3_b$-condition on the coefficient 
matrix $\sigma$.
(If $\sigma$ is not bounded, for example, 
then this parts becomes  difficult.)
Finally, concatenate these solutions on the subintervals
by using Chen's identity,
we obtain a time-global solution.

To prove local Lipschitz continuity, 
we estimate the distance between
$Z(m)$ and  $\hat{Z}(m)$ for $m \ge 1$ 
for two given rough paths  $X$ and $\hat{X}$ on each subinterval.
\qedsymbol

\bigskip

It is known that the $C_b^3$-condition on the coefficient 
matrix in Theorem \ref{thm_cont}
can be relaxed  to one called the
${\rm Lip} (\gamma)$-condition with $\gamma >p$
(see Lyons, Caruana and L\'evy \cite{lcl} for details).
Theorem \ref{thm_cont} also holds  for $p \ge 3$.
In that case, a suitable sufficient condition on the coefficient 
 is 
either $C^{[p]+1}_b$ or
${\rm Lip} (\gamma)$ with $\gamma >p$.

In fact, the Lyons-It\^o map is locally Lipschitz continuous
not just in $X$, 
but also in $\sigma$ and the initial condition.
So, it is quite flexible.
If particular, if  it is regarded as a map 
 in the initial condition (and the time) only,
it naturally defines a rough path version of 
a flow of diffeomorphism associated to an ODE/SDE.

There are other methods to solve RDEs.
Davie \cite{dav} solved an RDE 
by constructing a rough path version of 
the Euler(-Maruyama) approximation method
(see Friz and Victior \cite{fvbk} for details).
This method seems powerful.
Moreover, Bailleul  invented his flow method,
which is something like a "monster version" of 
Davie's method.
In this method, not just one initial value, but 
all initial values are considered simultaneously
and approximate solutions takes values in 
the space of homeomorphisms.
Using this, he recently solved RDEs with linearly growing coefficient.
(Precisely, the condition is that 
$\sigma$ itself may have linear growth,
but its derivatives are all bounded.
See \cite{be} for detail.)
Gubinelli's approach is also important,
which will be discussed in the next section.
However, it is not that his method to solve PDEs 
is different, 
but his formalism is.

%%%%%%%%%%%%%%%%%%%%%%%%%%%%%%%%%%%
\section{Gubinelli's controlled path theory}

The aim of this section is to present 
Gubinelli's formalism of rough path theory
in a nutshell.
It is recently called the controlled path theory
and currently competes with Lyons' original formalism.
It seems unlikely that one of the two
defeats the other in the near future.
However,
it is also true that this formalism is gaining attentions,
because it is simpler in a sense 
and  recently produced offsprings,
namely,
two new theories of singular stochastic PDEs
(Hairer's regularity structure theory
and Gubinelli-Imkeller-Perkovski's 
paracontrolled distribution theory).

The core of Gubinelli's idea is in his 
definition of rough path integrals. 
In Lyons' original definition, 
it is basically of the form $\int f(X)dX$.
In other words, 
"$X$ in $f(X)$"  and "$X$ in $dX$" must be the same
and 
an integral $\int f(Y)dX$ cannot be defined in general.
This is reasonable since 
a line integral is defined for a $1$-form.
However, 
impossibility of varying $X$ and $Y$ independently 
looks quite strange to
 most of probabilist who are not familiar 
with rough paths.
The author himself got surprised when 
he started studying the theory.
Since things like that are possible
in the Young (or Riemann-Stieltjes) integration 
and It\^o's stochastic integration,
some people may have tried  it for rough path integral 
only to find it hopeless.
Almost everyone gives up at this point 
and forget this issue.

Gubinelli did not, however.
He advanced "halfway"
by setting a Banach space of integrands 
for each rough path $X$ in an abstract way. 
Since this space contains elements of the form $f(X)$,
this is an extension of Lyons' rough path integration.
This Banach space depends on $X$
and may be different for different $X$.
Hence, this is not a complete separation 
of $X$ and the integrand.
Such an integrand is called a controlled path 
(with respect to $X$)
and $X$ is sometimes called a reference rough path.

One rough analogy for heuristic understanding 
is that
it looks like a "vector bundle" 
whose base space is an 
infinite dimensional curved space 
and whose fiber space is a Banach space.
The fiber spaces above different $X$'s are different 
vector spaces, although they look similar.

In this formalism, 
the integration map sends 
an integrand with respect to $X$ 
to an integral with respect to $X$
(which takes values in another Euclidean space)
for each fixed $X$.
Therefore, a solution to an RDE
driven by $X$ is understood as a fixed point 
in a certain Banach space of integrands with respect to $X$.

Now we give a brief mathematical explanation.
See Friz and Hairer \cite{fh} for details.
In this section we use $1/p$-H\"older topology
instead of  $p$-variation topology.
We assume $2 \le p <3$ for simplicity again, 
though the controlled path theory extends to 
the case $p \ge 3$.
The geometric rough path space 
with $1/p$-H\"older topology 
is denoted by $G\Omega^{H}_{1/p} ({\bf R}^d)$.
The $i$th level path of 
$X \in G\Omega^{H}_{1/p} ({\bf R}^d)$
is estimated by $i/p$-H\"older norm ($i=1,2$).

Let $X \in G\Omega^{H}_{1/p} ({\bf R}^d)$.
A pair $(Y, Y^{\prime})$ is said to be 
an ${\bf R}^n$-valued controlled path 
(controlled by $X$) 
if the following three conditions are satisfied:

\bigskip
\noindent
{\bf (i)}~$Y \in C^{1/p -Hld} ({\bf R}^n)$,
where
$C^{1/p -Hld} ({\bf R}^n)$ stands for the space of 
${\bf R}^n$-valued, $1/p$-H\"older continuous paths.
\\
{\bf (ii)}~$Y^{\prime} \in C^{1/p -Hld} ({\bf R}^{n } \otimes ({\bf R}^{d })^* )$.
\\
{\bf (iii)}~
If $R: \triangle \to {\bf R}^n$ is defined by 
$$
Y_{t}- Y_{s}= Y^{\prime}_s \cdot  X^1_{s,t} + R_{s,t}
\qquad
\qquad
(0 \le s \le t \le 1),
$$
then $R \in C^{2/p -Hld} (\triangle, {\bf R}^n)$ holds.

\bigskip

Note that ${}^{\prime}$ on the shoulder of $Y$
is just a symbol
and it does not mean differentiation with respect $t$.
Note also that $Y$ and $Y^{\prime}$ are 
one-parameter $1/p$-H\"older continuous functions,
while $R$ is a
two-parameter $2/p$-H\"older continuous function.
Loosely, the last condition means that
"behavior of $Y$ is at worst as bad as that of $X$"
since regularity of $R$ is better.
(To check this, fix $s$ arbitrarily and let $t$ vary near $s$.)
The totality of such $(Y, Y^{\prime})$ is 
denoted by ${\mathcal Q}_X^{1/p -Hld} ({\bf R}^n)$,
which becomes a Banach space 
equipped with the norm 
$\|Y\|_{1/p -Hld} + \|Y^{\prime}\|_{1/p -Hld} + \|R\|_{2/p -Hld}$.
One should note that
this Banach space of controlled paths 
depends on $X$.

Examples of controlled paths include:
{\bf (a)}~$X$ itself.
Precisely $t \mapsto X^1_{0,t}$.
{\bf (b)}~the composition $g(Y)$ for 
a $C^2$-function
$g: {\bf R}^{n} \to {\bf R}^{m}$
and $Y \in {\mathcal Q}_X^{1/p -Hld} ({\bf R}^{n})$.
{\bf (c)}~
The addition of 
$Y \in {\mathcal Q}_X^{1/p -Hld} ({\bf R}^{n})$
and a
$2/p$-H\"older continuous, ${\bf R}^n$-valued
path $Z$.
{\bf (d)}~
The multiplication of 
$Y \in {\mathcal Q}_X^{1/p -Hld} ({\bf R}^{n})$ and
a $2/p$-H\"older continuous, scalar-valued path $Z$.
(Precisely, the "derivatives" of these examples 
are naturally found and the pairs becomes controlled paths.)

For 
$(Y, Y^{\prime}) 
\in {\mathcal Q}_X^{1/p -Hld} ({\rm Mat}(n,d))$,
we can define a kind of 
the rough path integral along the reference rough path $X$
as the limit of a modified Riemann sum:
\[
(Z_{t} - Z_s:=) \quad
\int_s^t Y_u  dX_u 
= \lim_{|{\mathcal P}| \to 0} 
\sum_i
\Bigl\{
Y_{t_{i-1}} X^1_{t_{i-1}, t_i}  + Y^{\prime}_{t_{i-1}} \cdot X^2_{t_{i-1}, t_i}
\Bigr\}
\]
The second term of the summand on the right hand side 
is an element of ${\bf R}^n$
obtained as the contraction of 
$X^2_{t_{i-1}, t_i} \in ({\bf R}^d)^{\otimes 2}$
and 
$Y^{\prime}_{t_{i-1}} \in {\bf R}^n 
\otimes ({\bf R}^d)^*\otimes ({\bf R}^d)^*$.
If we set $Z^{\prime}_s= Y_s$, 
then we can show that
$(Z,  Z^{\prime}) \in {\mathcal Q}_X^{1/p -Hld} ({\bf R}^{n})$.
If $Y = f(X)$, $Z$ coincides with the first level path 
of the rough path integral in Lyons' original sense.
In this sense, 
rough path integration is generalized.

There is a significant difference, however.
In Lyons' formalism, 
the rough path integration maps 
a geometric rough path space to another.
In Gubinelli's formalism, 
it mops a controlled path space to another 
for each fixed reference geometric rough path $X$.
Moreover, it is linear in $(Y, Y^{\prime})$.
(One can prove continuity of 
integration map in $X$ in the latter formalism, too.)

Keeping these in mind, let us look at the 
driven ODE given at the beginning of this article:
\[
Y_t =  \int_0^t \sigma(Y_s) dX_s.
\]
Suppose that $Y \in {\mathcal Q}_X^{1/p -Hld} ({\bf R}^n)$. 
Then, 
composition $\sigma(Y) 
\in {\mathcal Q}_X^{1/p -Hld} ({\rm Mat}(n,d))$
and
the rough path integral belongs to 
$ {\mathcal Q}_X^{1/p -Hld} ({\bf R}^n)$ again.
Therefore, 
the integration map on the right hand side 
maps ${\mathcal Q}_X^{1/p -Hld} ({\bf R}^n)$
to itself
and it makes sense to think of its fixed points,
which are solutions to the RDE.
As we have seen, 
a solution to an RDE in this formalism is 
not a rough path, but a controlled path 
with respect to the driving rough path $X$.

If $\sigma$ is of $C_b^3$,
then the RDE has a unique global solution 
and the corresponding Lyons-It\^o map is 
locally Lipschitz continuous.
From this we can see that
the (first level paths of) the solutions 
in both Lyons' and Gubinelli's senses agree
for any $X$.
%

%%%%%%%%%%%%%%%%%%%%%%%
\section{Brownian rough path}

In the previous sections
everything was deterministic 
and no probability measure appeared so far.
In this section 
we lift the Wiener measure $\mu$
on the usual continuous path space
to 
a probability measure on the geometric rough path space,
by constructing a $G\Omega_p({\bf R}^d)$-valued
random variable called Brownian rough path.
In this theory Brownian rough path
plays the role of Brownian motion.

In this section
we assume $2<p<3$, excluding the case $p=2$.
We have denoted a (rough) path 
by $x$ or $X$ before.
To emphasize that it is a random variable
under the Wiener measure $\mu$, 
we will denote it by $w$ or $W$.

For 
$w \in C_0({\bf R}^d)$, 
denote by $w(m) \in C_0({\bf R}^d)$ the $m$th
dyadic polygonal approximation of $w$
associated with the partition $\{k/2^m~|~ 0 \le k \le 2^m\}$
($m=1,2,\ldots$).
Since it is clearly of finite variation, 
its natural lift $W (m)$ exists.
Set
$$
{\mathcal S} :=\{w \in C_0({\bf R}^d)~|~
\mbox{ $\{ W(m)\}_{m=1,2,\ldots}$ is  Cauchy 
in $G\Omega_p ({\bf R}^d)$}\}.
$$
Obviously,
a lift of $w \in {\mathcal S}$  is naturally defined naturally 
as
$\lim_{m \to \infty} W(m) \in G\Omega_p ({\bf R}^d)$.
If $w$
is of finite variation, 
then $w \in {\mathcal S}$ and the two kinds of lift 
actually agree.

How large is the subset ${\mathcal S}$?
In fact, it is of full Wiener measure. 
Hence, 
Brownian rough path $W$ can be  defined 
by this lift, that is,
$W:=\lim_{m \to \infty} W(m)$ $\mu$-a.s.
This is a $G\Omega_p ({\bf R}^d)$-valued random variable
defined on
$( C_0({\bf R}^d) , \mu)$
and its law (image measure) is a
probability measure on $G\Omega_p ({\bf R}^d)$.
(This lift map is neither deterministic nor continuous,
but is merely measurable.)
Thus, we  obtain
something like the Wiener measure on 
the geometric rough path space.
\footnote{In the author's view, this important measure deserves being named.}
By the way, this construction of Brownian rough path
works for $1/p$-H\"older topology, too.

Substituting $W$ into the Lyons-It\^o map,
we obtain a unique solution to 
the corresponding Stratonovich SDE.
Let us explain.
Consider RDE (\ref{eq.ode1})
and denote by
$\Phi: G\Omega_p ({\bf R}^d) \to G\Omega_p ({\bf R}^n)$
the associated  Lyons-It\^o map, namely,
 $Y =\Phi(W)$.

The SDE corresponding to (\ref{eq.ode1}) is given by
\begin{eqnarray}
dy_t  &=& \sigma (y_t) \circ dw_t
\nonumber\\
&=& \sigma (y_t)  dw_t+  
\frac12 {\rm Trace}[ \nabla \sigma (y_t)  \langle  \sigma (y_t) \bullet,  \bullet  \rangle ]dt, 
\qquad y_{0}=0.
%\label{eq.sde}
\nonumber
\end{eqnarray}
Compared to the SDE of It\^o-type, 
the above SDE has a modified term 
${\rm Trace} [\cdots]$ on the right hand side.
In terms of a Riemann sum, 
\[
\int_0^t \sigma (y_s) \circ dw_s = \lim_{|{\mathcal P}| \to 0} \sum_{i=1}^N   \frac{ \sigma (y_{t_{i}})+ \sigma( y_{t_{i-1}} ) } {2}  (  w_{t_{i}}  - w_{t_{i-1}}).
 \]
This is different from the Riemann sum 
for the corresponding It\^o integral.

\begin{theorem}\label{thm_sde}
Define  Brownian rough path $W$
as the lift of the canonical realization of 
Brownian motion $w =(w_t)_{0 \le t \le 1}$ as above.
Then, for almost all $w$ with respect to $\mu$,
$y_t = \Phi (W)^1_{0, t}$ holds for all $t \in [0,1]$.
\end{theorem}

This theorem states that 
a solution to an SDE  can be obtained
as the image of a continuous map.
It was inconceivable in usual probability theory.
The proof is easy.
Consider $\Phi( W(m))^1$ for each $m$. 
Since  RDEs 
are generalization of driven ODEs in 
the Riemann-Stieltjes sense, 
the unique solutions to 
the ODE driven by $w(m)$ and 
$\Phi( W(m))^1$ agree.
Then, take limits of both sides by using Wong-Zakai's
approximation theorem
and Lyons' continuity theorem,
which proves Theorem \ref{thm_sde}.

In the above argument,
the RDE and the corresponding 
SDE have no drift term, but
modification to the drift case is quite easy.
Instead of $W$, we just need to consider the Young pairing 
$(W, \lambda) \in G\Omega_p ({\bf R}^{d+1})$,
where 
$\lambda$ is the trivial one-dimensional path
given by $\lambda_t =t$.

In the end of this section
we discuss an application of Lyons' continuity theorem
to quasi-sure analysis 
(see Aida \cite{ai}, Inahama \cite{ina0, ina3}
for details).
Quasi-sure analysis is something like 
a potential theory on the Wiener space 
and one of the deepest topics in Malliavin calculus.
Those who are not familiar with Malliavin calculus 
can skip this part.

Since we have a suitable notion of differentiation 
on the Wiener space $(C_0({\bf R}^d) , \mu)$,
we can define a Sobolev space  ${\bf D}_{r,k}$,
where 
$r \in (1, \infty)$
 and 
 $k \in {\bf N}$
are the integrability and the differentiability indices,
respectively.
For each $(r, k)$ and subset $A$ of $C_0({\bf R}^d)$, 
we can define 
a capacity $C_{r,k} (A)$ 
via the corresponding Sobolev space.  
The capacity is finer than the Wiener measure $\mu$
and therefore a $\mu$-zero set may have positive capacity.
Since the
${\bf D}_{r,k}$-norm is increasing in both $r$ and $k$,
so is $C_{r,k} (A)$.
A subset $A$ is called slim if $C_{r,k} (A) =0$
for any $p$ and $k$.
Simply put, 
a slim set is much smaller than a typical  $\mu$-zero set.

Now we get back to rough path theory.
We have seen that ${\mathcal S}^c =0$,
but in fact we can prove that
${\mathcal S}^c$ is slim in a rather simple way.
Recall that $w \in {\mathcal S}^c$
is equivalent to that 
$w$ does not admit the lift via
the dyadic piecewise linear approximation.
Looking at the proof closely, 
we find that the lifting map $w \mapsto W$
is quasi-continuous.
(The results in this paragraph also hold
for $1/p$-H\"older topology.)

Consequently, 
the following famous theorems in
quasi-sure analysis becomes almost obvious.
First, the Wong-Zakai approximation 
theorem admits quasi-sure refinement.
This is now obvious since the lift map is defined outside a slim set.
Next, the solution to an SDE, 
as a path space-valued Wiener functional  
(or as a Wiener functional which takes values in 
the path space over flow of homeomorphism),
admits a
quasi-continuous modification.
This is also immediate from the quasi-continuity of the lift map
 (if we do not care about 
small difference of Banach norms).
Those who are not familiar with rough path theory
might be surprised that 
these results can be proved under the $C^3_b$-condition on $\sigma$,
not under smoothness.

%%%%%%%%%%%%%%%%%%%%%%%
\section{Gaussian rough path}\label{sec.gauss}

The aim of this section is 
to provide a summary of rough path lifts of 
Gaussian processes other than Brownian motion.
In this section roughness $p$ satisfies $2 \le p<4$.
This means that up to the third level paths, 
but not the fourth level path,
are to be considered. 
Such lifts of 
Gaussian processes are called Gaussian rough paths.
An RDE driven by a Gaussian rough path  is
something like an SDE driven by a Gaussian process.
Since an RDE is deterministic, 
whether 
the Gaussian process is a semimartingale or not 
is irrelevant.
If the Gaussian process admits a lift to a random rough path,
then we always have this kind of "SDE."

Let $w =(w^1_t, \ldots, w_t^d)_{0 \le t \le 1}$
be a $d$-dimensional, mean-zero 
Gaussian process with i.i.d. components.
We assume for simplicity that $w$ starts at the origin
so that it is a $C_0 ({\bf R}^d)$-valued random variable.
Its covariance is given by 
$R(s,t) :={\mathbb E}[w^1_s w^1_t]$
and determines the law of the process.

If $R(s,t)= \{ s^{2H} + t^{2H} - |t-s|^{2H} \}/2$
for some constant $H \in (0,1)$,
then $w$ is called fractional Brownian motion (fBm)
with the Hurst parameter $H$.
Generally, the smaller $H$ becomes, 
the tougher problems get.
When $H =1/2$, it is the usual Brownian motion.
When $H \neq 1/2$,
it is not a Markov process or a semimartingale anymore,
but it still has self-similarity and stationary increment.
From now on 
we only consider the case $1/4 < H \le 1/2$
unless otherwise stated.

A sample path of fBm is $1/p$-H\"older continuous
and of finite $p$-variation if $p >1/H$.
Hence, it is natural to ask 
whether $w$ admits a lift to a random rough path 
of roughness $p$ 
as in the previous section.
By using the dyadic polygonal approximations,
Coutin and Qian \cite{cq} lifted fBm $w$
to a $G\Omega_p ({\bf R}^d)$-valued random variable $W$
if $1/4 < H\le 1/2$ and $p>1/H$.
(The smaller $p$ is, the stronger 
the statement becomes.
So, this condition should be understood as 
"for $p$ slightly larger than $1/H$.")
When $1/3 <H \le 1/2$,
we can take $[p]=2$ and only use 
the first and the second level paths.
When $1/4 <H \le 1/3$, however,
$[p]=3$ and the third level path is needed.
The lift $W$ is called fractional Brownian rough path.
It is the first Gaussian rough path discovered 
and is still the most important example
(other than Brownian rough path).
This result also holds for $1/p$-H\"older topology.
By the way, 
this kind of rough path lift fails when $H \le 1/4$.
(Even in such a case, a "non-standard" lift 
of fBm exists.)

Can we lift  more general Gaussian processes 
beyond special examples such as fBm?
Since the covariance $R(s,t)$ knows everything about
the Gaussian process $w$,
it seems good to impose certain conditions on $R(s,t)$.
However, it is not easy to find a suitable 
sufficient condition.
Friz and Victoir \cite{fv} noticed that 
$\rho$-variation norm $\|R\|_{\rho}$ 
of $R$ as a two-parameter 
function should be considered.
We set for $\rho \ge 1$, 
\[
\|R\|_{\rho}^{\rho} 
=
\sup_{{\mathcal P}, {\mathcal Q}} \sum_{i, j}  \Bigl|
R(s_{i},  t_{j})
- R(s_{i},  t_{j-1})
- R(s_{i-1},  t_{j})
+R(s_{i-1},  t_{j-1})
\Bigr|^{\rho}.
\]
Here, the supremum runs over all the pairs of
partitions of 
${\mathcal P}=(s_i)$ and ${\mathcal Q}=(t_j)$ of $[0,1]$.

Let us consider the natural lift $\{W(m)\}_{m=1,2,\ldots}$ of 
the dyadic piecewise linear approximations of $w$
as in the previous section.
According to \cite{fv}, 
if $1 \le \rho <2$ and $2 \rho< p<4$,
then each level path
$\{W(m)^i\}$ converges in $L^r$ for every 
$r \in (1, \infty)$ as a sequence of random variables
which take values in the Banach space of 
$p/i$-variation topology ($1 \le i \le [p]$).
The limit $W$ is call Gaussian rough path
or a lift of $w$.
If $R$ satisfies some kind of H\"older condition
in addition (as in the case of fBm),
then convergence takes place in 
$i/p$-H\"older topology,
and, moreover, convergence is not just $L^r$, but 
also almost sure.

We have only discussed the lift via the dyadic piecewise 
linear approximations.
However, it is proved that many kinds of lift 
in fact coincide.
Examples include the mollifier method, 
a more general piecewise linear approximations,
the Karhunen-Lo\'eve approximations
(approximation of $w$ by a linear combination 
of an orthonormal basis of Cameron-Martin space
and i.i.d. of one-dimensional standard normal distributions).
In this sense, $W$ is a canonical lift of $w$,
though not unique.

Thus, Gaussian rough path $W$ exists if  $1 \le \rho <2$.
The next question is how nice $W$ is.
If the lift map destroys structures of the Gaussian measure, 
then studying solutions to RDEs driven by $W$ 
may become very difficult.
Among many structures on a Gaussian space, 
the most important one is probably 
Cameron-Martin theorem.  
It states that the image measure of the Gaussian measure
induced by a translation along Cameron-Martin vector
is mutually absolutely continuous to  the Gaussian measure.
Therefore, one  naturally hopes that 
the rough path lifting should not destroy 
the structure of translations. 
This is what the complementary Young regularity condition is about.
Loosely, it demands that 
the Cameron-Martin translation on the lower space 
(i.e., the abstract Wiener space)
and the Young translation on the upper space 
(i.e., the geometric rough path space)
should be compatible.

As above, suppose that 
$R$ is of finite $\rho$-variation 
for some $\rho \in [1,2)$.
Cameron-Martin space of the Gaussian process $w$
is denoted by ${\mathcal H}$.
We say that the complementary Young regularity is satisfied
if there exist $p$ and $q$ with the following properties:
$p \in (2\rho, 4)$, $q \in [1,2)$, $1/p +1/q >1$
and ${\mathcal H}$ is continuously embedded in
$C_0^{q -var} ({\bf R}^d)$,
the set of continuous paths 
of finite $q$-variation starting at $0$.

In this case, since $W={\mathcal L}(w)$ takes values in 
$G\Omega_p  ({\bf R}^d)$, 
the Young translation by 
an element of $C_0^{q -var} ({\bf R}^d)$ is well-defined.
Here, ${\mathcal L}$ denotes the rough path lift map.
If the complementary Young regularity holds,
then there exists a subset $A$ of full measure 
such that for any $h \in {\mathcal H}$ and $w \in A$,
${\mathcal L}(w+h) = \tau_h ({\mathcal L}(w))$ holds.
In other words, 
lifting and translation commute.
Thanks to this nice property, 
we can prove many theorems under  
the complementary Young regularity condition,
as we will see.

In the case of Brownian motion, 
Cameron-Martin paths behave nicely.
However, it is not easy to study 
behavior of 
Cameron-Martin paths for the other Gaussian processes
including fBm.
Friz and coauthors \cite{fvimb, fggr} proved that, 
for fBm with Hurst parameter $H\in (1/4, 1/2]$,  
${\mathcal H}\subset
C_0^{q -var} ({\bf R}^d)$ with $q =(H +1/2)^{-1}$.
Since we can take any $p >1/H$,
we can find $p$ and $q$ with $1/p +1/q >1$,
which is the condition for Young integration.
Therefore, fBm satisfies 
the complementary Young regularity condition
if $H\in (1/4, 1/2]$.

Other
examples of sufficient condition 
for  the complementary Young regularity is as follows.
(1)~$R$ is of finite $\rho$-variation
for some $\rho \in [1, 3/2)$.
(2)~A quantity called the mixed $(1,\rho)$-variation
of $R$ is finite for $\rho \in [1,2)$ (see \cite{fggr} for the latter).

In the author's opinion, 
the simplest way to understand
the current theory of Gaussian rough paths is
as follows
(a class with a larger number is smaller):
(i)~The covariance $R$ is of finite $\rho$-variation
for some $\rho \in [1, 2)$.
In this case, the canonical rough path lift exists.
(ii)~The case that the complementary 
Young regularity condition is  satisfied in addition.
(iii)~fBm with $H \in (1/4,1/2]$ as the most important example.

%%%%%%%%%%%%%%%%%%%%%%%%%%%%%%%%%%%%%%%%%%%%%%%%%%%%%%%%
\section{Large deviation principle}
\label{sec.ldp}

From now on 
we will review probabilistic results in rough path theory.
The aim of this section is to discuss
a Schilder-type
large deviation principle (LDP).

Let us recall the standard version of 
Schilder's LDP for Brownian motion.
Let $\mu$ be the Wiener on $C_0({\bf R}^d)$
and let ${\mathcal H}$ be Cameron-Martin space.
We denote by $\mu^{\ve}$ 
the image measure of the scalar multiplication map 
$w \mapsto \ve w$ by $\ve >0$.
A good rate function $I: C_0({\bf R}^d) \to [0, \infty]$
is defined by $I(w) = \|w\|^2_{{\mathcal H}}/2$
for $w \in {\mathcal H}$
and $I(w) =\infty$ for $w \notin {\mathcal H}$.
Obviously, the mass concentrates at the origin
as $\ve \searrow 0$.
Moreover, the following Schilder's LDP holds:
\[
- \inf_{w \in A^{\circ}}  I(w)
\le
\liminf_{\ve \searrow 0}\ve^2 \log \mu^{\ve} (A^{\circ}) \le  \limsup_{\ve \searrow 0} \ve^2\log \mu^{\ve} (\bar{A})
\le 
- \inf_{w \in \bar{A}} I(w)
\]
for every Borel subset $A \subset C_0({\bf R}^d)$,
where $A^{\circ}$ and $\bar{A}$
denote the interior and the closure of $A$, respectively.
Roughly, this claims that 
weight of the subset $A$ which is distant from $0$
decays like $\exp(- \mbox{const}/\ve^2)$
and the positive constant can be written as 
the infimum of the rate function $I$ over $A$.
A little bit mysterious is that 
information on ${\mathcal H}$ dictates the LDP 
though $\mu^{\ve}({\mathcal H})=0$ for any $\ve >0$.

LDPs go well with continuous maps.
If an LDP holds on a domain of a continuous map, 
then it is transferred to an LDP on the image 
and the new rate function can be written 
in terms of the original one.
This is called the contraction principle.
Let us take a look at Freidlin-Wentzell's LDP
from this viewpoint.

For a sufficiently nice coefficient matrix 
$\sigma: {\bf R}^n \to {\rm Mat}(n,d)$
and a drift vector $b: {\bf R}^n \to {\bf R}^n$,
consider the following Stratonovich-type SDE
index by $\ve >0$:
\[
dy^{\ve}_t = \sigma(y_t^{\ve}) \circ \ve dw_t +b (y^{\ve}_t)dt   
\qquad \mbox{with \quad $y^{\ve}_0 =0\in {\bf R}^n$.}
\]
One can easily guess that 
the law of the stochastic process $y^{\ve}$
concentrates around a unique solution to 
the following deterministic ODE 
"$dz_t= b (z_t)dt$ with $z_0=0$."
In fact, a stronger result, Freidlin-Wentzell's LDP holds.

Formally, 
if $\Phi$ denotes the {\it usual} It\^o map 
associated with the block matrix 
$[\sigma, b]$ and the initial value $y_0 =0$
and 
$\lambda_t =t$, 
then $y^{\ve} = \Phi (\ve w, \lambda)$.
Therefore, 
if the usual It\^o map {\it were} continuous,
Freidlin-Wentzell's LDP would be 
immediate from Schilder's LDP and the contraction principle.
In reality, $\Phi$ is not continuous.
Hence, this LDP was proved by other methods.
\footnote{A few methods are known.}
However, the rigorously proved statement 
is the same as the one obtained by 
the  formal argument as above.

In such a situation, Ledoux, Qian and Zhang \cite{lcl}
gave a new proof of this LDP 
using rough path theory.
Let $\hat\mu^{\ve}$ be the law of the scaled 
Brownian rough path $\ve W=(\ve W^1, \ve^2 W^2)$.
First, they proved a Schilder-type LDP for $\{\hat\mu^{\ve} \}$
on $G\Omega_p ({\bf R}^d)$ for $2<p<3$.
More precisely, 
for every Borel subset $A \subset  G\Omega_p ({\bf R}^d)$,
\[
- \inf_{X \in A^{\circ}}  \hat{I}(X)
\le
\liminf_{\ve \searrow 0}\ve^2 \log \hat\mu^{\ve} (A^{\circ}) 
\le  \limsup_{\ve \searrow 0} \ve^2\log \hat\mu^{\ve} (\bar{A})
\le 
- \inf_{X \in \bar{A}} \hat{I}(X)
\]
holds.
Here, $\hat{I}$ is a good rate function on 
$G\Omega_p ({\bf R}^d)$ defined by 
$\hat{I}(X) = \|h\|^2_{{\mathcal H}}/2$ 
if $X = {\mathcal L}(h)$ for some $h \in {\mathcal H}$
and 
$\hat{I}(X) =\infty$ if otherwise
(${\mathcal L}$ is the rough path lift).

Next, recall that we rigorously have
$y^{\ve} = [ t \mapsto \Phi (\ve W, \lambda)^1_{0,t}]$
in rough path theory.
Here, $(\ve W, \lambda)$ is the Young pairing 
and 
$\Phi^1$ is the first level of 
the Lyons-It\^o map, which is continuous.
Hence, 
we can actually use the contraction principle 
to prove Freidlin-Wentzell's LDP.

This work attracted attention 
because of its clear perspective on the LDP
and its making use of Lyons' continuity theorem.
Many papers followed it 
and there are now many variants 
of Schilder-type LDPs on the geometric rough path space.
A prominent example is one for Gaussian rough paths.
If $\|R\|_{\rho} <\infty$ for some $1 \le \rho <2$,
then a Schilder-type LDP holds 
for the laws of the scaled Gaussian rough path 
on $G\Omega_p ({\bf R}^d)$ with $p> 2\rho$
(see Theorem 15.55, \cite{fvbk}).
Of course, the case of 
$1/p$-H\"older topology is studied as well.
Now it seems that
the Schilder-type LDP on rough path space
has become an independent topic itself,
separating from the original motivation
of showing Freidlin-Wentzell's LDP.

%%%%%%%%%%%%%%%%%%%%%%%%%%%%%%%%%%%%%%%%%%%%%%%%%%%%%%
\section{Support theorem}

Like Freidlin-Wentzell's LDP, 
one could easily prove Stroock-Varadhan's 
support theorem if the usual It\^o map were continuous.
The aim of this section is to summarize 
Ledoux-Qian-Zhang's a new proof in \cite{lcl}
via rough path.

We consider the following SDE with the same 
coefficients $\sigma$ and $b$:
\[
dy_t = \sigma(y_t) \circ  dw_t +b (y_t)dt   
\qquad \mbox{with \quad $y_0=0 \in {\bf R}^n$.}
\]
The solution $y=(y_t)_{0 \le t \le 1}$
induces an image measure on
$C_{0}({\bf R}^n)$.
What is its support, that is, the smallest closed 
subset which carries the whole weight?

The support theorem answers this question.
It claims that we should look at where 
the corresponding deterministic It\^o map
sends Cameron-Martin paths.
For $h \in {\mathcal H}$, let $\phi (h)$ be a unique 
solution to the following driven ODE: 
\begin{equation}\label{ode.phi}
d\phi (h)_t = \sigma(\phi(h)_t) dh_t +b (\phi(h)_t)dt   
\qquad \mbox{with \quad $\phi(h)_0=0 \in {\bf R}^n$.}
\end{equation}
Then, the support is the closure of 
$\{ \phi(h) | h \in {\mathcal H} \}$ in $C_{0}({\bf R}^n)$.
The support of the Wiener measure $\mu$ is 
the domain of It\^o map $C_{0}({\bf R}^d)$
and 
${\mathcal H}$ is  dense in it.
Hence, if It\^o map {\it were} continuous, 
the support theorem would be very easy.
However, the proof was hard in reality.

Since Lyons-It\^o map is continuous
and  extends the deterministic It\^o map,  
the support theorem is immediate 
if one checks the support of the law of 
Brownian rough path $W$.
In fact, they proved that the support is 
the closure of ${\mathcal L} ({\mathcal H})$
(the lift of $\mathcal H$ in $G\Omega_p ({\bf R}^d)$),
which is actually the whole set $G\Omega_p ({\bf R}^d)$ ($2<p<3$).
From this
Stroock-Varadhan's support theorem follows at once.

The support theorem on the geometric 
rough path space was generalized 
to the case of Gaussian rough paths
with complementary Young regularity condition
(see Theorem 15.60, \cite{fvbk}).
The case of $1/p$-H\"older topology was also studied.

%%%%%%%%%%%%%%%%%%%%%%%%%%%%%%%%%%%%%%%%%%%%%%%%%%%%%%
\section{Laplace approximation}

In this section 
we discuss the Laplace approximation,
\footnote{There could be small differences among the literature
in what the terms like 
Laplace approximation, Laplace asymptotics, Laplace's method
precisely mean.}
 that is, 
the precise asymptotics  of the LDP 
of Freidlin-Wentzell type in Section \ref{sec.ldp}.
We consider the case where 
the driving rough path is fractional Brownian rough path.
(See \cite{ina2}. 
The case of infinite dimensional Brownian rough path 
is in \cite{ik}.)

Consider the same RDE as in Section \ref{sec.ldp}:
\[
dy^{\ve}_t = \sigma(y_t^{\ve})  \ve dx_t +b (y^{\ve}_t)dt   
\qquad \mbox{with \quad $y^{\ve}_0 =0\in {\bf R}^n$,}
\]
where $\ve \in (0,1]$ is a small parameter.
%
%For simplicity, the initial value is $0$.
%
As a driving rough path $X$, 
we take fractional Brownian rough path $W$
with Hurst parameter $H\in (1/4,1/2]$.
Take $p \in (1/H, [1/H] +1)$.
From the Schilder-type LDP for the law of $\ve W$
on $G\Omega_p ({\bf R}^d)$ 
and Lyons' continuity theorem, 
the law of the solution 
$y^{\ve} =[t \mapsto (Y^{\ve})^1_{0, t}]$
also satisfies LDP 
of Freidlin-Wentzell type on $C_0^{p -var} ({\bf R}^n)$.

Let
$\phi(h)$ be the solution to ODE (\ref{ode.phi})
for $h \in {\mathcal H}$, but
${\mathcal H} ={\mathcal H}^H$
stands for Cameron-Martin space for fBm here.
As we have seen, $h \in {\mathcal H}$ 
is of finite $q$-variation with
$q = (H+1/2)^{-1}~(<2)$.
Hence, this driven ODE should be understood in the Young sense.

By a general fact called Varadhan's lemma, 
the following limit theorem holds:
\[
\lim_{\ve \searrow 0}\ve^2 
\log {\mathbb E}
\bigl[ \exp \bigl(- F (y^{\ve}) /\ve^2 \bigr) 
\bigr]
=
- \inf_{h \in {\mathcal H}} 
\bigl\{
F(\phi(h)) + \frac12 \|h  \|_{{\mathcal H}}^2 
\bigr\}
\]
for every bounded continuous function
$F: C_0^{p-var} ({\bf R}^n) \to {\bf R}$.
This is an "integral form" of the LDP of Freidlin-Wentzell type.

The above formula calculates the logarithm 
of a certain expectation of exponential type.
The Laplace approximation studies asymptotic 
behavior of the expectation of exponential type itself 
under additional assumptions on $F$.
The case of the usual SDE was first proved 
by Azencott \cite{az} and Ben Arous \cite{ba},
followed by many others.
\footnote{The strongest ones among them 
show a kind of 
Laplace approximation in the framework on Malliavin calculus
to obtain asymptotics of heat kernels.}

In this article, we consider this problem 
from a viewpoint of rough path theory.
An advantage of this approach is as follows.
The most important part of the proof is Taylor expansion 
of (Lyons-)It\^o map.
It also becomes deterministic in rough path theory
and therefore small difference of the original Gaussian process
does not matter as long as it admits a rough path lift.
Consequently, the rough path proof 
can treat the cases of the usual Brownian motion
and fBm in a unified way.

Now we introduce assumptions:

\vspace{3mm}
\noindent
{\bf (H1):}
For some $p >1/H$, 
$F$ and $G$ are real-valued
bounded continuous functions defined on
$C^{p -var}_0 ({\bf R}^n)$.

\vspace{3mm}
\noindent
{\bf (H2):}
A real-valued 
function $\hat{F} $ on ${\mathcal H}$ defined by 
$\hat{F} := F \circ \phi   +  \|  \,\cdot\, \|^2_{{\mathcal H} } /2$
achieves a minimum exactly at one point
$\gamma \in {\mathcal H}$.

\vspace{3mm}
\noindent
{\bf (H3):}
On a certain neighborhood of $\phi(\gamma)$ in 
$C^{p -var}_0 ({\bf R}^n)$,
$F$ and $G$ are Fr\'echet smooth and 
all of their derivatives are bounded.

\vspace{3mm}
\noindent
{\bf (H4):}
The Hessian 
$\nabla^2 (F \circ \phi ) (\gamma)   
\vert_{{\mathcal H} \times {\mathcal H}}$
of $F \circ \phi \vert_{\mathcal H}$
at $\gamma \in {\mathcal H}$
is strictly larger than 
$-\langle \,\cdot\, , \,\cdot\ \rangle_{{\mathcal H}}$ 
in the form sense.

\vspace{3mm}

These assumptions are typical 
for Laplace approximations.
We assume in addition that the coefficients 
$\sigma$ and $b$ are bounded, 
smooth with bounded derivatives of all order.
Then, we can show the following 
asymptotic expansion:

As $\ve \searrow 0$ we have
\begin{eqnarray*}
{\mathbb E}  \bigl[  
G(  y^{\ve})   \exp \bigl(   - F (  y^{\ve})   /\ve^2  \bigr)
\bigr]
=
\exp ( - \hat{F} (\gamma)  /\ve^2 ) 
\bigl(
\alpha_0  +\alpha_1 \ve + \cdots +   \alpha_m \ve^m  + \cdots
\bigr)
\end{eqnarray*} 
for certain constants $\alpha_j~(j=0,1,2,\ldots)$.

A key of proof is a Taylor-like expansion of 
(the first level of)
the Lyons-It\^o map on a neighborhood of
the lift of $\gamma$ in $G\Omega_p ({\bf R}^d)$.
This expansion is deterministic 
and irrelevant to any probability measure
or stochastic process.
By the way, we can see from the LDP that 
contributions from the complement set 
of the neighborhood is negligible.

A more detailed explanation is as follows.
Denote by $\Phi:G\Omega_p ({\bf R}^{d+1}) \to G\Omega_p ({\bf R}^n)$
the Lyons-It\^o map associated with the 
coefficient $[\sigma | b]$.
Substitute the Young pairing
$(\ve X,\lambda) \in G\Omega_p ({\bf R}^{d+1})$
 of 
$\ve X \in G\Omega_p ({\bf R}^{d})$
and $\lambda_t =t$ into $\Phi$.
Then, $y^{\ve}_t =\Phi((\ve X,\lambda))^1_{0,t}$ and
$\phi(h)_t =\Phi((h,\lambda))^1_{0,t}$,
where $h$ and its natural lift is denoted by the same symbol.
What we need in the proof is an
expansion of $\Phi((\gamma +\ve X,\lambda))^1$.
Note that $\gamma +\ve X$ is 
in fact a Young translation by $\gamma$.

There exist $\phi_j(\gamma, X)~(j=1,2,\ldots)$
such that the first level path admits the following 
expansion as $\ve \searrow 0$
with respect to the $p$-variation topology:
\[
\Phi((\gamma +\ve X,\lambda))^1
=
\phi(\gamma)
+
\ve \phi_1 (\gamma, X)+ \ve^2 \phi_2 (\gamma, X) +\cdots.
\]
Formally, each $\phi_j(\gamma, X)$ satisfies
a simple ODE of first order 
and can be written down by the variation of constants formula.
Since $\phi_j(\gamma, X)$
is of order $j$ as a functional of $X$,
the above expansion is something like Taylor expansion.
Of course, the remainder term
also satisfies a reasonable estimate
with respect to the rough path topology.
(One should note here that, 
while $X$ can be  an arbitrary element in 
$G\Omega_p ({\bf R}^d)$,
$\gamma$ has to be a "nice" path so that 
the Young translation works.
Otherwise, this expansion would not make sense.)

Let us take a look from a slightly different angle.
This explanation might be easier 
for non-experts of rough path theory.
In the setting of finite variation, 
the It\^o map in the Riemann-Stieltjes sense 
is known to be Fr\'echet smooth.
Hence, it admits a Taylor expansion 
around any $\gamma$ 
for an infinitesimal vector $\ve w$.
The Taylor terms are formally the same 
as $\phi_j$'s above.
The Taylor-like expansion for $\Phi^1$
is a completion with respect to the rough path topology
 of this Taylor expansion in Fr\'echet sense.

This Taylor-like expansion also holds when $p \ge 3$.
In that case, the base point $\gamma$ 
can be of $q$-variation with any $q\in [1,2)$
such that $1/p +1/q >1$.
(See \cite{ina1}).

%%%%%%%%%%%%%%%%%%%%%%%%%%%%%%%%%%%%%%%%%%%%%%%
\section{Jacobian processes and their moments}

As before we consider the SDE 
with the coefficients $\sigma$ and $b$,
but we denote 
the column vectors of $\sigma$ 
by $V_1, \ldots, V_d$
and $b$ by $V_0$.
We should regard them as
vector fields on ${\bf R}^n$.
In this section $V_i~(0 \le i \le d)$ are assumed to be of $C^{[p]+2}_b$.
Using this notation, 
we can rewrite the RDE with a general initial condition
as
\begin{equation}\label{wai}
dy_t = \sum_{i=1}^d V_i (y_t)  dx^i_t + V_0 (y_t)dt   
\qquad \mbox{with \quad $y_0=a \in {\bf R}^n$.}
\end{equation}
Here,  the superscript $i$ on the shoulder of 
$dx^i$ stands for the coordinate of ${\bf R}^n$,
not the level of an iterated integral.

Take (formal) differentiation of 
$y_t=y_t(a)$ with respect to the initial value $a \in {\bf R}^n$.
Then, $j_t = \nabla y_t $ and $k_t =j_t^{-1}$
are $n \times n$-matrices and satisfy the 
following ODEs at least formally:
\begin{eqnarray}
d j_t &=&  \sum_{i=1}^d \nabla V_i (y_t) j_t dx^i_t + \nabla V_0 (y_t) j_t dt, 
\qquad
\mbox{with \quad $j_0 ={\rm Id}_n$.}
\label{jay}
\\
d k_t &=&  - \sum_{i=1}^d k_t \nabla V_i (y_t)  dx^i_t - k_t \nabla V_0 (y_t)  dt, 
\qquad
\mbox{with \quad $k_0 ={\rm Id}_n$.}
\label{kay}
\end{eqnarray}
Here, $\nabla V_i$ is regarded as a $n \times n$-matrix, too.
This $j$ is called a Jacobian process of 
the original differential equation 
and plays a very important role in analysis 
of the usual SDEs.
Therefore, it should be very important 
in rough path theory, too.

When we regard 
 (\ref{wai})-(\ref{kay}) as a system of RDEs 
 driven by a geometric rough path $X$,
we have troubles.
The first one is that the coefficients of (\ref{jay})-(\ref{kay})
are not bounded.
Hence, we cannot use the standard version 
of Lyons' continuity theorem (Theorem \ref{thm_cont}).
RDEs with unbounded coefficients 
are often difficult to handle and their solutions
may explode in finite time.

However, 
the system of RDEs (\ref{wai})-(\ref{kay})
has a unique time-global solution
and Lyons' continuity theorem holds.
The reason is as follows.
Because of a "triangular" structure of the system, 
RDE (\ref{wai}) solves first and we obtain $(X,Y)$.
Now that $\nabla V_i (y_t)$ is  known,
(\ref{jay}) and (\ref{kay}) become linear RDEs.
A solution to a linear driven ODE 
can be expressed as an infinite sum. 
By generalizing this argument to the case 
of rough path topology,
we obtain the first level path of 
a solution to a linear RDE.
If the first level path of the solution stays
inside a sufficiently large ball, 
then behavior of the coefficients 
outside the ball is irrelevant
and we can use a standard cut-off technique
to obtain higher level paths of the solution
and prove the continuity theorem.
Thus, we have seen 
 that for every $X \in G\Omega_p ({\bf R}^d)$,
the system (\ref{wai})-(\ref{kay}) has a unique 
time-global solution $(Y,J,K)$.

In stochastic analysis for RDEs, 
integrability of $J$ and $K$ matters.
Due to the cut-off argument as above, 
it is sufficient to prove integrability of 
$\sup_{0 \le t \le 1} (| J^1_{0, t}| + | K^1_{0, t}| )$.
However, this was quite hard since 
a straight forward computation yields 
\[
\sup_{0 \le t \le 1} (| J^1_{0, t}| + | K^1_{0, t}| )
\le C \exp \Bigl( C \sum_{i=1}^{[p]}  \|X^i\|_{p/i}^{p/i} \Bigr)
\]
for some constant $C>0$.
If $X$ is a Gaussian rough path, then we usually have $p>2$.
Fernique's theorem is not available
and the right hand side is not even in $L^1$.

An integrability lemma by 
Cass, Litterer and Lyons \cite{cll} solves this problem.
For any $\alpha >0$, set $\tau_0 =0$
and 
\[
\tau_m = 1 \wedge \inf\{ t \ge \tau_{m-1} ~|~ \sum_{i=1}^{[p]}  \|X^i\|_{p/i, [ \tau_{m-1} ,t]}^{p/i}
\ge \alpha \}
\]
recursively for $m \ge 1$.
Here, 
$\|X^i\|_{p/i, [s,t]}$ is the $p/i$-variation norm of $X^i$
restricted on the subinterval $[s,t]$.
Then, we set 
\[
N_{\alpha} (X) = \max\{m~|~ \tau_m <1 \}.
\]
This quantity is important.
Since it is non-increasing in $\alpha$,
integrability  of $N_{\alpha}$ is valuable
for small $\alpha$.

If we compute on each subinterval $[\tau_{m-1}, \tau_{m}]$, 
we can prove 
\begin{equation}
\label{estJK}
\sup_{0 \le t \le 1} (| J^1_{0, t}| + | K^1_{0, t}| )
\le C_{\alpha} \exp \Bigl( C_{\alpha} N_{\alpha} (X) \Bigr),
\end{equation}
where $C_{\alpha}>0$ is a constant which may 
depend on $\alpha$.
Note that this is a deterministic estimate.
Therefore, it is sufficient to show the 
exponential integrability of $N_{\alpha}$ for some $\alpha$.

They proved in \cite{cll} that 
for a Gaussian rough path  $W$
with the complementary Young regularity condition,
there exists $\delta >0$ such that
${\mathbb E}[ \exp(  N_{\alpha}(W)^{1+\delta} ) ] <\infty$
for any $\alpha >0$.
(More precisely, they gave a sharp estimate 
of the tail probability of $N_{\alpha}(W)$).
Consequently, both sides of 
(\ref{estJK}) have moments of all order.
Thus, we have obtained
moment estimates of the Jacobian process 
for an RDE driven by a Gaussian rough path.

%%%%%%%%%%%%%%%%%%%%%%%%%%%%%%%%%%%%%%%%%%%%%%%%%
\section{Malliavin calculus for rough differential equations}

As most of successful theories in analysis,
Malliavin calculus has its abstract part 
and concrete examples of functionals
to which the abstract theory apply.
The former is the theory of Sobolev spaces
on an abstract Wiener spaces.
In other words, it is differential and integral calculus
on an infinite dimensional Gaussian space.
The latter is  solutions to SDEs.
Hence, 
a natural question is whether Malliavin calculus 
is applicable to RDEs driven by a Gaussian rough path.  
In this section, we consider RDE (\ref{wai})
with $C_b^{\infty}$-coefficient vector fields
driven by  a Gaussian rough path as in 
Section \ref{sec.gauss}.

The first works in this direction are 
Cass and Friz (and Victoir) \cite{cfv, cf}.
They consider the lift of Gaussian processes 
with the complementary Young condition 
and a certain non-degeneracy condition.
(Loosely, 
this "non-degeneracy" condition above 
is to exclude Gaussian processes
that do not diffuse very much
such as the pinned Brownian motion.)
Examples include fractional Brownian motion 
with $H \in (1/4, 1/2]$.
They proved that the solution $y_t$ to the RDE
is differentiable in a weak sense, 
namely, 
it belongs to the local Sobolev space 
${\bf D}_{r, 1}^{loc}({\bf R}^n)~(1<r <\infty)$.
Moreover, if $V_i~(0 \le i \le d)$ satisfies 
H\"ormander's bracket generating condition 
at the starting point $y_0 =a\in {\bf R}^n$,
Malliavin covariance matrix of $y_t$
is non-degenerate in a weak sense, 
namely, 
it is invertible a.s.,
which implies that the law of $y_t$
has a density $p_t(a,a')$ with respect to
the Lesbegue measure $da'$.
However, this argument bring us no information on
regularity of the density.

As in the study of the usual SDEs, 
we would like to show that {\rm (i)}
the solution $y_t$ belongs to the Sobolev space
${\bf D}_{r, k}({\bf R}^n)$
for any integrability index
$r \in (1,\infty)$ and the differentiability index $k\ge 0$
and
{\rm (ii)} Malliavin covariance matrix of $y_t$
is non-degenerate in the sense of Malliavin, 
namely, 
the determinant of the inverse 
of the covariance matrix  has moments of all order. 
These imply smoothness of the density in $a'$. 
The biggest obstacle was 
the moment estimates of the Jacobian process
as in the previous section, however.
For example, 
$D^k y_t$,
the $k$th derivative of the solution,
 has an explicit expression which involves 
 the Jacobian processes and its inverse.
Therefore, unless this obstacle was removed, 
we could not proceed.

After the moment estimates was recently proved,
Malliavin calculus for RDEs has developed rapidly.
First, Hairer and Pillai \cite{hp} proved the case of fBm 
with $H \in (1/3, 1/2]$.
Differentiability in the sense of Malliavin calculus,
that is, 
$y_t \in 
{\bf D}_{r, k}({\bf R}^n)$
for any $k \ge 0$ and $1 <r <\infty$,
was shown by fractional calculus.
The point in their proof of 
non-degeneracy of the Malliavin covariance matrix 
is a deterministic version of Norris' lemma
in the framework of the controlled path theory.

Under the Young complementary regularity condition,
differentiability was shown in Inahama \cite{ina4}.
The theory of Wiener chaos,
not fractional calculus, is used in the proof.
Non-degeneracy under H\"ormander's condition 
was proved by 
Cass, Hairer, Litterer and Tindel \cite{chlt}
(and Baudoin, Ouyang and Zhang \cite{boz1})
for a rather general class of Gaussian processes
including fBm with $H \in (1/4, 1/2]$.
Since these recent results enables us to 
study RDEs with Malliavin calculus quite smoothly,
this research topic may make great 
advances  in the near future.

Some papers on this topic in the case of fBm 
were already published.
(1)~Varadhan's estimate, that is, 
short time asymptotics of 
the logarithm of the density $\log p_t(a,a')$
(see \cite{boz1}).
(2)~Smoothing property of the "heat semigroup"
under Kusuoka's UFG condition 
(see Baudoin, Ouyang and Zhang \cite{boz2}).
This condition is on the Lie brackets of 
the coefficient vector fields and weaker than 
H\"ormander's condition.
(3)~Positivity of the density 
$p_t(a,a')$ 
(see Baudoin, Nualart, Ouyang and Tindel \cite{bnot}).
In these three papers, $1/4 < H\le 1/2$,
while in the next paper $1/3 < H\le 1/2$.
(4)~Short time off-diagonal asymptotic expansion 
of the density $p_t(a,a')$
under the ellipticity assumption 
on the coefficients at the starting point
(see Inahama \cite{ina5}).
In the last example, 
Watanabe's theory of generalized Wiener functionals,
(that is, Watanabe distributions) 
and asymptotic theory for them are used.
The theory is known to be a very powerful tool 
in Malliavin calculus,
but it also works well in 
the frameworks of rough path theory.
We also point out that the proof of the off-diagonal asymptotics 
is a kind of Laplace approximation in the framework of Malliavin calculus
and therefore the Taylor-like expansion in Section 11 plays a crucial role.

%%%%%%%%%%%%%%%%%%%%%%%%%%%
\section{Topics that were not covered}
For lack of space, 
we did not discuss some important topics
in and around rough path theory.
The most important among them is 
applying ideas from rough path theory 
to stochastic partial differential equations 
(SPDEs).
This is an attempt to use rough path theory
to solve 
 singular SPDEs which cannot be solved 
by existing methods.
(One should not misunderstand that 
the general theory of SPDE is rewritten 
or extended with rough paths.)
Several attempts have already been published,
but there seems to be no unified theory.
So, reviewing them in details
in a short article like this is impossible,
but we give a quick comment on two of them  which look very active now.

The two most successful ones are 
Hairer's regularity structure theory \cite{ha2, fh}
and Gubinelli-Imkeller-Perkowski's
para-controlled distribution thoery \cite{gip, gp}.
\footnote{Another example is "fully nonlinear rough stochastic PDEs" 
studied by P. Friz and coauthors.}
Examples of singular SPDEs 
these theories solve are similar, including
KPZ equation, the dynamic $\Phi^4_3$,
three dimensional stochastic Navier-Stokes equation, etc.,
but the two theories look quite different.
They now should be classified as independent theories, 
not as a part of rough path theory.

The numerical and the statistical studies of the usual SDEs 
are very important.
Hence, it might be interesting to consider 
analogous problems for RDEs driven by
 a Gaussian rough path.
Not so many papers have been written by now, 
but we believe these topics will be much larger.
Approximations of SDEs from  a viewpoint 
of rough paths should be included in this paragraph, too.

Neither have we mentioned signatures of rough path
T. Lyons and coauthors study intensively.
For a (rough) path defined on the time interval $[0,1]$,
its iterated integrals on the whole interval 
$$
X^k_{0,1}=\int_{ 0 < t_1 <\cdots <t_k <1} dx_{t_1}\otimes \cdots \otimes dx_{t_k}
\qquad
(k=1,2,\ldots)
$$
(or  the corresponding quantities) are called signature 
of the (rough) path.
A fundamental problem in this topic is 
whether the signatures determine a (rough) path 
modulo reparametrization.
A probabilistic version is whether 
the expectations of the signatures determine
a probability measure on the (rough) path space.
For recent results, see \cite{ly2} and references therein.


\begin{thebibliography}{00}


\bibitem{ai}
Aida, S.;  Vanishing of one-dimensional $L^2$-cohomologies of loop groups. 
J. Funct. Anal. 261 (2011), no. 8, 2164--2213.

\bibitem{az}
Azencott, R.;
Formule de Taylor stochastique et d\'eveloppement asymptotique d'int\'egrales de Feynman. 
 Seminar on Probability, XVI, Supplement,  pp. 237--285, Lecture Notes in Math., 921, Springer, Berlin-New York, 1982.



\bibitem{be}
Bailleul, I.;
Flows driven by rough paths.
 Rev. Mat. Iberoam. 31 (2015), no. 3, 901--934.

\bibitem{boz1}
Baudoin, F.; Ouyang, C.;  Zhang, X.;
Varadhan Estimates for rough differential equations driven by fractional Brownian motions.
 Stochastic Process. Appl. 125 (2015), no. 2, 634--652. 


\bibitem{boz2}
Baudoin, F.; Ouyang, C.;  Zhang, X.;
Smoothing effect of rough differential equations driven by fractional Brownian motions.
Ann. Inst. Henri Poincar\'e Probab. Stat. 52 (2016), no. 1, 412--428.


\bibitem{bnot}
Baudoin, F.; Nualart, E.; Ouyang, C.; Tindel. S.;
On probability laws of solutions to differential systems driven by a fractional Brownian motion.
To appear in Ann. Probab.
 arXiv:1401.3583


\bibitem{ba}
Ben Arous, G.; 
Methods de Laplace et de la phase stationnaire sur l'espace de Wiener.   Stochastics  25  (1988),  no. 3, 125--153.



\bibitem{cfo} 
Caruana, M.; Friz, P.; Oberhauser, H.;
A (rough) pathwise approach to a class of non-linear stochastic partial differential equations. 
Ann. Inst. H. Poincar\'e Anal. Non Lin\'eaire 28 (2011), no. 1, 27--46. 



\bibitem{cf}
 Cass, T.; Friz, P.; 
 Densities for rough differential equations under H\"ormander's condition. 
 Ann. of Math. (2) 171 (2010), no. 3, 2115--2141. 

\bibitem{cfv}
 Cass, T.; Friz, P.; Victoir, N.; 
 Non-degeneracy of Wiener functionals arising from rough differential equations. 
 Trans. Amer. Math. Soc. 361 (2009), no. 6, 3359--3371.



\bibitem{chlt}
Cass, T.; Hairer, M.; Litterer, C.;  Tindel, S.;
Smoothness of the density for solutions to Gaussian rough differential equations. 
Ann. Probab. 43 (2015), no. 1, 188--239. 


\bibitem{cll}
Cass, T.; Litterer, C.; Lyons, T.;
 Integrability and tail estimates for Gaussian rough differential equations. 
 Ann. Probab. 41 (2013), no. 4, 3026--3050. 

\bibitem{cq}
Coutin, L.; Qian, Z.;
Stochastic analysis, rough path analysis and fractional Brownian motions.  
Probab. Theory Related Fields  122  (2002),  no. 1, 108--140.




\bibitem{dav}
 Davie, A. M.; 
 Differential equations driven by rough paths: an approach via discrete approximation. 
 Appl. Math. Res. Express. AMRX 2007, no. 2, Art. ID abm009, 40 pp.





\bibitem{fggr}
Friz. P.; Gess, B.; Gulisashvili, A.; Riedel, S.;
Jain-Monrad criterion for rough paths and applications.
To appear in Ann. Probab.
  arXiv:1307.3460.

\bibitem{fh} 
 Friz, P.; Hairer, M.;
A course on rough paths. With an introduction to regularity structures. 
  Springer,  2014.


\bibitem{fvimb}
 Friz, P.; Victoir, N.;
A variation embedding theorem and applications. 
J. Funct. Anal. 239 (2006), no. 2, 631--637.
 
\bibitem{fv} 
 Friz, P.; Victoir, N.; 
 Differential equations driven by Gaussian signals. 
 Ann. Inst. Henri Poincar\'e Probab. Stat. 46 (2010), no. 2, 369--413.




\bibitem{fvbk}
Friz, P.; Victoir, N.;
Multidimensional stochastic processes as rough paths.
 Cambridge University Press, Cambridge, 2010. 

\bibitem{gu}
Gubinelli, M. 
Controlling rough paths. 
J. Funct. Anal. 216 (2004), no. 1, 86--140. 


\bibitem{gip}
Gubinelli, M.; Imkeller, P.; Perkowski, N.;
%
Paracontrolled distributions and singular PDEs.
Forum Math. Pi 3 (2015), e6, 75 pp.
 
 
\bibitem{gp}
Gubinelli, M.;  Perkowski, N.;
KPZ reloaded.
Preprint. arXiv:1508.03877 
%
 
\bibitem{ha1}
 Hairer, M.; 
 Solving the KPZ equation. 
 Ann. of Math. (2) 178 (2013), no. 2, 559--664.


\bibitem{ha2}
 Hairer, M.; 
A theory of regularity structures. 
 Invent. Math. 198 (2014), no. 2, 269--504. 


\bibitem{hp}
 Hairer, M.; Pillai, N.; 
 Regularity of laws and ergodicity of hypoelliptic SDEs driven by rough paths. 
 Ann. Probab. 41 (2013), no. 4, 2544--2598

%\bibitem{hl}
% Hambly, B.; Lyons, T.; 
% Uniqueness for the signature of a path of bounded variation and the reduced path group. 
% Ann. of Math. (2) 171 (2010), no. 1, 109--167.




\bibitem{ina0}
Inahama, Y.;
Quasi-sure existence of Brownian rough paths and a construction of Brownian pants.
Infin. Dimens. Anal. Quantum Probab. Relat. Top. 9 (2006), No. 4, 513--528. 

\bibitem{ina1}
Inahama, Y.;
A stochastic Taylor-like expansion in the rough path theory,
J. Theoret. Probab. 23 (2010), Issue 3, 671--714. 

\bibitem{ina2} 
Inahama, Y.;
Laplace approximation for rough differential equation driven by fractional Brownian motion.
Ann. Probab. 41 (2013), No. 1, 170-205. 

\bibitem{ina3} 
Inahama, Y.;
Large deviation principle of Freidlin-Wentzell type for pinned diffusion processes.
Trans. Amer. Math. Soc. 367 (2015), no. 11, 8107--8137.

\bibitem{ina4}
Inahama, Y.;
Malliavin differentiability of solutions of rough differential equations.
J. Funct. Anal. 267 (2014), no. 5, 1566--1584.


\bibitem{ina5}
Inahama, Y.;
Short time kernel asymptotics for rough differential equation driven by fractional Brownian motion.
Preprint. arXiv:1403.3181.


\bibitem{ik}
Inahama, Y.; Kawabi, H.;
Asymptotic expansions for the Laplace approximations for It\^o functionals of Brownian rough paths.
J. Funct. Anal. 243 (2007), no. 1, 270--322. 


\bibitem{lqz}
Ledoux, M.; Qian, Z.; Zhang, T.; 
Large deviations and support theorem for diffusion processes via rough paths. 
Stochastic Process. Appl. 102 (2002), no. 2, 265--283.

% \bibitem{lejq}
% Le Jan, Y.; Qian, Z.;
%Stratonovich's signatures of Brownian motion determine Brownian sample paths. 
%Probab. Theory Related Fields 157 (2013), no. 1-2, 209--223.  
 
\bibitem{ly}
Lyons, T.; 
Differential equations driven by rough signals. 
Rev. Mat. Iberoamericana 14 (1998), no. 2, 215--310. 


\bibitem{ly2}
Lyons, T.; 
Rough paths, signatures and the modelling of functions on streams.
To appear in the Proceedings of the International Congress of Mathematicians 2014, Korea.
arXiv: 1405.4537.




\bibitem{lcl}
 Lyons, T.; Caruana, M.; L\'evy, T.; 
Differential equations driven by rough paths. 
  Lecture Notes in Math., 1908. Springer, Berlin, 2007.

\bibitem{lq}
 Lyons, T.; Qian, Z.; 
 System control and rough paths. 
 Oxford University Press, Oxford, 2002. 


\bibitem{su}
Sugita, H.; 
%Various topologies in the Wiener space and L\'evy's stochastic area. 
%Probab. Theory Related Fields 91 (1992), no. 3-4, 283--296.
%
Hu-Meyer's multiple Stratonovich integral and essential continuity of multiple Wiener integral. 
Bull. Sci. Math. 113 (1989), no. 4, 463--474. 

 \end{thebibliography}
\end{document}